\documentclass[12pt]{article}
\usepackage{latexsym,amssymb,amsmath,geometry,setspace,enumerate}
\newtheorem{theorem}{Theorem}
\newtheorem{conjecture}[theorem]{Conjecture}
\newtheorem{lemma}[theorem]{Lemma}
\newtheorem{corollary}[theorem]{Corollary}
\newtheorem{proposition}[theorem]{Proposition}

\usepackage{lineno}

\begin{document}
\onehalfspace

\title{\Large Local Connectivity, Local Degree Conditions, some Forbidden Induced Subgraphs, and Cycle Extendability}
\author{Christoph Brause$^1$, Dieter Rautenbach$^2$, and Ingo Schiermeyer$^1$}
\date{}
\maketitle
\vspace{-10mm}
\begin{center}
{\small 
$^1$ Institute of Discrete Mathematics and Algebra, TU Bergakademie Freiberg\\ 
Freiberg, Germany, $\{$brause, schiermeyer$\}$@math.tu-freiberg.de\\[3mm]
$^2$ Institute of Optimization and Operations Research, Ulm University\\
Ulm, Germany, dieter.rautenbach@uni-ulm.de}
\end{center} 

\begin{abstract}
The research in the present paper was motivated by the conjecture of Ryj\'{a}\v{c}ek 
that every locally connected graph is weakly pancyclic.

For a connected locally connected graph $G$ of order at least $3$, 
our results are as follows:
If $G$ is $(K_1+(K_1\cup K_2))$-free, then $G$ is weakly pancyclic.
If $G$ is $(K_1+(K_1\cup K_2))$-free, then $G$ is fully cycle extendable if and only if $2\delta(G)\geq n(G)$.
If $G$ is 
$\{ K_1+K_1+\bar{K}_3,K_1+P_4\}$-free
or 
$\{ K_1+K_1+\bar{K}_3,K_1+(K_1\cup P_3)\}$-free,
then $G$ is fully cycle extendable.
If $G$ is distinct from $K_1+K_1+\bar{K}_3$ 
and
$\{ K_1+P_4,K_{1,4},K_2+(K_1\cup K_2)\}$-free,
then $G$ is fully cycle extendable.

Furthermore, if $G$ is a connected graph of order at least $3$ such that 
$$|N_G(u)\cap N_G(v)\cap N_G(w)|>|N_G(u)\setminus (N_G[v]\cup N_G[w])|$$
for every induced path $vuw$ of order $3$ in $G$, 
then $G$ is fully cycle extendable,
which implies that every connected 
locally Ore or locally Dirac graph of order at least $3$ is fully cycle extendable.
\end{abstract}

{\small 

\noindent \begin{tabular}{lp{12.5cm}}
\textbf{Keywords:} & 
Local connectivity; 
forbidden induced subgraph; 
cycle; 
weakly pancyclic; 
fully cycle extendable\\
\textbf{MSC2010:} &
05C38; 05C40; 05C45
\end{tabular}

}

\pagebreak

\section{Introduction}

We consider finite, simple, and undirected graphs, and use standard terminology.
A graph $G$ is {\it locally connected} if for every vertex $u$ of $G$,
the subgraph $G[N_G(u)]$ of $G$ induced by the neighborhood $N_G(u)$ of $u$ in $G$ is connected.
Similarly, $G$ is {\it locally Ore} if 
$|N_G(u)\cap N_G(v)|+|N_G(u)\cap N_G(w)|\geq d_G(u)$ for every induced path $vuw$ of order $3$ in $G$,
and $G$ is {\it locally Dirac} if 
$2\delta(G[N_G(u)])\geq d_G(u)$ for every vertex $u$ of $G$,
where $\delta(G)$ and $d_G(u)$ are the minimum degree of $G$ and the degree of $u$ in $G$, respectively.
For a vertex $u$ of a graph $G$ and some positive integer $k$, 
let $N_G^k(u)$ be set of vertices of $G$ at distance exactly $k$ from $u$.
For a graph $G$,
the {\it girth} $g(G)$ and the {\it circumference} $c(G)$
are the minimum and the maximum order of a cycle in $G$,
respectively.
A graph $G$ is {\it hamiltonian} if $c(G)=n(G)$,
where $n(G)$ is the order of $G$.
A graph $G$ is {\it weakly pancyclic} if it has a cycle of order $\ell$ for every integer $\ell$ between $g(G)$ and $c(G)$.
A cycle $C$ in a graph $G$ is {\it extendable}
if $G$ contains a $C'$ of order $n(C)+1$ with $V(C)\subseteq V(C')$.
A graph $G$ is {\it fully cycle extendable} if every vertex of $G$ lies on a triangle,
and every cycle in $G$ of order less than $n(G)$ is extendable.
For two disjoint graphs $G$ and $H$, 
let $G\cup H$ be their union,
let $G+H$ be their join, and
let $\bar{G}$ be the complement of $G$.
Let $K_n$, $P_n$, and $C_n$ be the complete graph, path, and cycle of order $n$, respectively.
Note that the star $K_{1,n-1}$ of order $n$ is $K_1+\bar{K}_{n-1}$.
Let $K_{1,1,n-2}=K_1+K_1+\bar{K}_{n-2}$.
If ${\cal F}$ is a set of graphs, then a graph is {\it ${\cal F}$-free}
if it does not contain a graph in ${\cal F}$ as an induced subgraph.
If ${\cal F}$ contains only one graph $F$, 
we write {\it $F$-free} instead of ${\cal F}$-free. 

The results in the present paper are motivated by the following conjecture.

\begin{conjecture}[Ryj\'{a}\v{c}ek \cite{r}]\label{conj_1}
Every locally connected graph is weakly pancyclic.
\end{conjecture}
We identify several sets ${\cal F}$ of graphs 
such that Conjecture \ref{conj_1} holds for ${\cal F}$-free graphs.
Furthermore, we show that Conjecture \ref{conj_1} holds for graphs 
that are locally Ore or locally Dirac.

Before we proceed to our results, 
we give a very short summary of previous related work.
Chartrand and Pippert \cite{cp} showed that every connected locally connected graph $G$
with $n(G)\geq 3$ and maximum degree $\Delta(G)\leq 4$ is either hamiltonian or $K_{1,1,3}$.
Extending a result of Kikust \cite{k},
Hendry \cite{h} showed that every connected locally connected graph $G$ 
with $n(G)\geq 3$, $\Delta(G)\leq 5$, and $\Delta(G)-\delta(G)\leq 1$ is fully cycle extendable.
Gordon et al. \cite{gops} generalized this last result further to graphs $G$ 
with $n(G)\geq 3$, $\Delta(G)\leq 5$, and $\delta(G)\geq 3$.
Extending earlier results due to Oberly, Sumner \cite{os}, and Clark \cite{c},
Zhang \cite{z} showed that every connected locally connected $K_{1,3}$-free graph $G$
with $n(G)\geq 3$ is fully cycle extendable.
Faudree et al. \cite{frs} weakened the local connectivity requirement for this last result.
Kubicka et al. \cite{kubicka} considered locally Dirac graphs,
and showed that every connected locally Dirac graph $G$ 
with $n(G)\geq 3$ and $\Delta(G)\leq 11$ is fully cycle extendable.
Furthermore, they observed that results of Hasratian and Khachatrian \cite{hk} 
imply that every connected locally Dirac graph of order at least $3$ is hamiltonian.
Our results generalize the mentioned results of 
Zhang \cite{z}
and 
Kubicka et al. \cite{kubicka}.

\section{Results}

Our first goal is to verify Conjecture \ref{conj_1} for $(K_1+(K_1\cup K_2))$-free graphs.
The next lemma collects several useful observations.

\begin{lemma}\label{lemma1}
Let $G$ be a connected locally connected $(K_1+(K_1\cup K_2))$-free graph of order at least $3$.
\begin{enumerate}[(i)]
\item Every vertex of $G$ lies on a triangle.
\item The diameter of $G$ is at most $2$.
\item $N_G^2(u)$ is independent for every vertex $u$ of $G$.
\item A cycle $C$ in $G$ is extendable if and only if there is some vertex $u$ in $V(G)\setminus V(C)$
that has a neighbor in $V(C)$ such that $N_G(u)\not\subseteq V(C)$ or $2d_G(u)>n(C)$.
\item Every cycle in $G$ of length less than $c(G)$ is extendable.
\end{enumerate}
\end{lemma}
{\it Proof:} (i) Since $G$ is connected locally connected and distinct from $K_2$,
the minimum degree of $G$ is at least $2$, which implies that every vertex of $G$ lies on a triangle.

\medskip

\noindent (ii) For a contradiction, we assume that $uvxy$ is some shortest path in $G$.
Since the minimum degree of $G$ is at least $2$, and $G$ is locally connected,
the vertices $u$ and $v$ have a common neighbor $w$.
Since $G$ is $(K_1+(K_1\cup K_2))$-free, $x$ is adjacent to $w$.
Since $y$ is not adjacent to $v$ or $w$, 
the graph $G[\{ v,w,x,y\}]$ is $K_1+(K_1\cup K_2)$,
which is a contradiction. Hence, $G$ has diameter at most $2$.

\medskip

\noindent (iii) For a contradiction, we assume that $xy$ be an edge between two vertices in $N_G^2(u)$
for some vertex $u$ of $G$.
Since $G$ is $(K_1+(K_1\cup K_2))$-free,
$x$ and $y$ have no common neighbor in $N_G(u)$.
Let $v$ be a common neighbor of $u$ and $x$,
and
let $w$ be a common neighbor of $u$ and $y$.
Note that $x$ is not adjacent to $w$,
and that $y$ is not adjacent to $v$.
Since $G$ is $(K_1+(K_1\cup K_2))$-free,
$v$ and $w$ are not adjacent.
Since $G$ is locally connected,
$v$ and $u$ have a common neighbor $z$.
Since $G$ is $(K_1+(K_1\cup K_2))$-free,
we obtain that 
$x$ and $z$ are adjacent,
and that
$y$ and $z$ are adjacent,
that is, $z$ is a common neighbor of $x$ and $y$ in $N_G(u)$,
which is a contradiction.
Hence, $N^2_G(u)$ is independent for every vertex $u$ of $G$.

\medskip

\noindent (iv) Let $C$ be a cycle in $G$.

First, we prove the necessity.
For a contradiction, we may assume that $C$ is extendable but that 
$N_G(u)\subseteq V(C)$ and $2d_G(u)\leq n(C)$ for every vertex $u$ in $V(G)\setminus V(C)$.
Since $C$ is extendable, there is some vertex $u$ in $V(G)\setminus V(C)$ 
such that $G$ contains a cycle $C'$ with $V(C')=V(C)\cup \{ u\}$.
By (ii) and (iii), at least one of every two consecutive vertices of $C$ is adjacent to $u$.
Since $2d_G(u)\leq n(C)$, this implies that exactly one of every two consecutive vertices of $C$ is adjacent to $u$,
that is, $C$ is a cycle of even length that alternates between $N_G(u)$ and the independent set $N_G^2(u)$,
where $|N_G(u)|=|N_G^2(u)|$.
Since $C'$ contains two edges between $u$ and $N_G(u)$,
the independence of $N_G^2(u)$ implies the contradiction $|N_G(u)|=|N_G^2(u)|+1$.

Next, we prove the sufficiency.
Therefore, let $u$ in $V(G)\setminus V(C)$
be such that $u$ has a neighbor in $V(C)$ and $N_G(u)\not\subseteq V(C)$ or $2d_G(u)>n(C)$.
If $N_G(u)\subseteq V(C)$, 
then $2d_G(u)>n(C)$ implies that $u$ is adjacent to two consecutive vertices of $C$, 
say $v$ and $w$, 
and replacing $vw$ with $vuw$ yields a cycle $C'$ 
of order $n(C)+1$ with $V(C)\subseteq V(C')$.
Hence, we may assume that $N_G(u)\not\subseteq V(C)$.
Since $u$ has a neighbor on $C$ and a neighbor not on $C$, 
the local connectivity of $G$ implies that there is a triangle $uvwu$
with $v\in V(G)\setminus V(C)$ and $w\in V(C)$.
Since $G$ is $(K_1+(K_1\cup K_2))$-free,
$u$ or $v$ is adjacent to a neighbor of $w$ on $C$,
and we obtain a cycle $C'$ similarly as above.

\medskip

\noindent (v) For a contradiction, we assume that $C$ is a cycle in $G$ of length less than $c(G)$
such that $C$ is not extendable.
By (iv), $V(G)\setminus V(C)$ is an independent set, 
and $2d_G(u)\leq n(C)$ for every vertex $u$ in $V(G)\setminus V(C)$.
Let $u^*\in V(G)\setminus V(C)$.
By (iii), at least one of every two consecutive vertices of $C$ is adjacent to $u^*$.
Since $2d_G(u^*)\leq n(C)$, this implies that exactly one of every two consecutive vertices of $C$ is adjacent to $u^*$,
that is, $C$ is a cycle of even length that alternates between $N_G(u^*)$ and the independent set $N_G^2(u^*)$,
where $|N_G(u^*)|=|N_G^2(u^*)|$.
By (ii) and symmetry, we obtain that $N_G(u)=N_G(u^*)$ for every vertex $u$ in $V(G)\setminus V(C)$,
that is, $V(G)\setminus N_G(u^*)$ is an independent set.
This implies that $c(G)\leq 2|N_G(u^*)|=n(C)$,
which is a contradiction.
$\Box$

\medskip

\noindent With Lemma \ref{lemma1} at hand, 
it is easy to verify Conjecture \ref{conj_1} for $(K_1+(K_1\cup K_2))$-free graphs.

\begin{theorem}\label{theorempaw}
Let $G$ be a connected locally connected $(K_1+(K_1\cup K_2))$-free graph of order at least $3$.
\begin{enumerate}[(i)]
\item $G$ is weakly pancyclic.
\item $G$ is fully cycle extendable if and only if $2\delta(G)\geq n(G)$.
\end{enumerate}
\end{theorem}
{\it Proof:} (i) follows immediately from Lemma \ref{lemma1} (i) and (v).
We proceed to the proof of (ii).
If $2\delta(G)\geq n(G)$, then the theorem of Dirac \cite{d} implies $c(G)=n(G)$,
and (ii) follows from Lemma \ref{lemma1} (i) and (v).
Now, let $G$ be fully cycle extendable.
Let $u$ be a vertex of $G$ of minimum degree.
By Lemma \ref{lemma1} (ii) and (iii), $V(G)\setminus N_G(u)$ is an independent set.
Since $G$ has a hamiltonian cycle, 
we obtain $n(G)-d_G(u)=|V(G)\setminus N_G(u)|\leq |N_G(u)|=d_G(u)$,
which implies $2d_G(u)\geq n(G)$.
$\Box$

\medskip

\noindent We proceed to further connected locally connected graphs that are fully cycle extendable.

In our next result we consider forbidding just one induced subgraph.
Let $X$ be the graph with vertex set $\{ u_0,u_1,u_2,u_3,u_1',u_2',u_3'\}$ 
and edge set 
$$\{ u_0u_1,u_0u_2,u_0u_3,u_0u_1',u_0u_2',u_0u_3',u_1u_2,u_2u_3,u_1u_3,u_1u_1',u_2u_2',u_3u_3'\}.$$

\begin{proposition}\label{proposition1}
Let $F$ be a graph.

Every connected locally connected $F$-free graph of order at least $3$ is fully cycle extendable
if and only if $F$ is an induced subgraph of $K_{1,3}$ or $K_1+P_3$.
\end{proposition}
{\it Proof:} Let $G$ be a connected locally connected $F$-free graph of order at least $3$.

First, we prove the sufficiency.
If $F$ is an induced subgraph of $K_{1,3}$, then Zhang's result \cite{z} implies that $G$ is fully cycle extendable.
If $F$ is an induced subgraph of $K_1+P_3$, 
then $N_G(u)$ induces a complete graph for every vertex $u$ of $G$,
which implies that $G$ is complete, and hence $G$ is fully cycle extendable.

Next, we prove the necessity.
Since $K_{1,1,3}$ and $X$ are connected locally connected graphs of order at least $3$ that are not hamiltonian,
$F$ must be an induced subgraph of $K_{1,1,3}$ as well as of $X$.
Since $K_{1,1,3}$ is not an induced subgraph of $X$, 
$F$ is a proper induced subgraph of $K_{1,1,3}$,
which implies that $F$ is an induced subgraph of $K_{1,3}$ or $K_1+P_3$.
$\Box$

\medskip

\noindent The arguments used in the previous proof lead to the following result 
concerning pairs of forbidden induced subgraphs.

\begin{proposition}\label{proposition2}
Let ${\cal F}$ be a set of two graphs that contains no induced subgraph of $K_{1,3}$ or $K_1+P_3$.

If every connected locally connected ${\cal F}$-free graph of order at least $3$ is fully cycle extendable,
then one graph in ${\cal F}$ is $K_{1,1,3}$, 
and the other graph in ${\cal F}$ is an induced subgraph of $X$.
\end{proposition}
{\it Proof:} Since $K_{1,1,3}$ is not hamiltonian, 
and ${\cal F}$ contains no proper induced subgraph of $K_{1,1,3}$,
the set ${\cal F}$ must contain $K_{1,1,3}$.
Since $X$ is $K_{1,1,3}$-free and not hamiltonian, 
the other graph in ${\cal F}$ is an induced subgraph of $X$.
$\Box$

\medskip

\noindent The next two results yield examples for sets ${\cal F}$ as in Proposition \ref{proposition2}.
Theorem \ref{theorem1} can actually be derived from 
Theorem \ref{theorem6} below, 
but we include a short independent proof 
using the forbidden induced subgraphs.
Note that 
Theorem \ref{theorem4}
extends the result of Zhang \cite{z}, 
because both, 
$K_1 + (K_1 \cup P_3)$ and $K_{1,1,3}$,
contain $K_{1,3}$ as an induced subgraph. 

\begin{theorem}\label{theorem1}
Every connected locally connected $\{ K_{1,1,3},K_1+P_4\}$-free graph of order at least $3$ 
is fully cycle extendable.
\end{theorem}
{\it Proof:} 
Let $G$ be a connected locally connected $\{ K_{1,1,3},K_1+P_4\}$-free graph of order at least $3$.
As before, 
the minimum degree of $G$ is at least $2$, and every vertex of $G$ lies on a triangle.
For a contradiction, we may assume that $C$ is a cycle in $G$ of order less than $n(G)$
such that $G$ contains no cycle $C'$ of order $n(C)+1$ with $V(C)\subseteq V(C')$.
We fix a cyclic order on $C$.
For every vertex $u$ on $C$, 
let $u^-$ and $u^+$ be the predecessor and successor of $u$ on $C$ within the cyclic order.
Since $G$ is connected, some vertex $u$ on $C$ has a neighbor $x$ in $V(G)\setminus V(C)$.
Our assumption implies that $x$ is adjacent to neither $u^-$ nor $u^+$.
Since $G[N_G(u)]$ is connected and $P_4$-free, $G[N_G(u)]$ contains a path $xvu^-$.
Our assumption implies that $v$ lies on $C$.

First, we assume that $u^-$ is adjacent to $u^+$ for every vertex $u$ that has a neighbor in $V(G)\setminus V(C)$.
Since $v$ is adjacent to $x$, $v^-$ is adjacent to $v^+$.
Note that $v^+=u^-$ or $v^-=u^+$ is possible.
Since $G[\{ u,u^-,u^+,x,v\}]$ is not $K_1+P_4$, $v$ is adjacent to $u^+$,
and the cycle $uu^+\ldots v^-v^+\ldots u^-vxu$ contradicts our assumption.
Hence, we may assume that $u^-$ is not adjacent to $u^+$.

Since $G[\{ u,u^-,u^+,x,v\}]$ is not $K_{1,1,3}$, $v$ is not adjacent to $u^+$.
Since $G[N_G(u)]$ is connected and $P_4$-free, $G[N_G(u)]$ contains a path $xwu^+$
where $w$ is distinct from $v$.
By symmetry of $v$ and $w$, $w$ is not adjacent to $u^-$.
If $v$ and $w$ are adjacent, then $G[\{ u,u^-,u^+,v,w\}]$ is $K_1+P_4$, and,
if $v$ and $w$ are not adjacent, then $G[\{ u,u^-,x,v,w\}]$ is $K_1+P_4$, 
which is a contradiction, and completes the proof. $\Box$

\begin{theorem}\label{theorem4}
Every connected locally connected $\{ K_{1,1,3},K_1+(K_1\cup P_3)\}$-free graph of order at least $3$ 
is fully cycle extendable.
\end{theorem}
{\it Proof:} 
Let $G$ be a connected locally connected $\{ K_{1,1,3},K_1+(K_1\cup P_3)\}$-free graph of order at least $3$.
As before,
the minimum degree of $G$ is at least $2$, and every vertex of $G$ lies on a triangle.
For a contradiction, we may assume that $C$ is a cycle in $G$ of order less than $n(G)$
such that $G$ contains no cycle $C'$ of order $n(C)+1$ with $V(C)\subseteq V(C')$.
We fix a cyclic order on $C$.
For every vertex $u$ on $C$, 
let $u^-$ and $u^+$ be the predecessor and successor of $u$ on $C$ within the cyclic order.
Since $G$ is connected, some vertex $u$ on $C$ has a neighbor $x$ in $V(G)\setminus V(C)$.
Our assumption implies that $x$ is adjacent to neither $u^-$ nor $u^+$.
Let $d$ be the minimum distance within the graph $G[N_G(u)]$ 
between the vertex $x$ and a vertex in $\{ u^-,u^+\}$.
Clearly, $d\geq 2$.
Since $G[N_G(u)]$ is connected and $(K_1\cup P_3)$-free, we have $d\leq 3$.

First, we assume that $u^-$ and $u^+$ are not adjacent.
If $d=2$, then, 
by symmetry, we may assume that $u$, $u^-$, and $x$ have a common neighbor $v$.
Now, 
if $v$ is not a neighbor of $u^+$, then $G[\{ u,u^-,u^+,v,x\}]$ is $(K_1+(K_1\cup P_3))$, 
and,
if $v$ is a neighbor of $u^+$, then $G[\{ u,u^-,u^+,v,x\}]$ is $K_{1,1,3}$,
which is a contradiction.
Hence, $d=3$, and, by symmetry, we may assume that $G[N_G(u)]$ contains a path $xvwu^-$.
Now, 
if $w$ is not a neighbor of $u^+$, then $G[\{ u,u^+,v,w,x\}]$ is $(K_1+(K_1\cup P_3))$,
and,
if $w$ is a neighbor of $u^+$, then $G[\{ u,u^-,u^+,v,w\}]$ is $K_{1,1,3}$,
which is a contradiction.
Hence, $u^-$ and $u^+$ are adjacent.

Next, we assume that $d=2$.
By symmetry, we may assume that $u$, $u^-$, and $x$ have a common neighbor $v$.
By our assumption, the vertex $v$ lies on $C$. 
By symmetry between $u$ and $v$, $v^-$ and $v^+$ are adjacent.
Note that $v^-=u^+$ or $v^+=u^-$ is possible.
Now the cycle $uu^+\ldots v^-v^+\ldots u^-vxu$ contradicts our assumption.
Hence, $d=3$.

By symmetry, we may assume that $xvwu^-$ is a path in $G[N_G(u)]$.
By our assumption, the vertex $w$ lies on $C$. 
Clearly, $v$ is adjacent to neither $u^-$ nor $u^+$.
If $v$ does not lie on $C$, 
then the minimum distance 
between the vertex $v$ and a vertex in $\{ u^-,u^+\}$ is $2$, 
in which case we can argue as above for $d=2$.
Hence, we may assume that $v$ lies on $C$.
By symmetry, $v^-$ and $v^+$ are adjacent.
If $w$ is not a neighbor of $u^+$, then $G[\{ u,u^-,u^+,w,x\}]$ is $(K_1+(K_1\cup P_3))$,
which is a contradiction.
Hence, $w$ is a neighbor of $u^+$.
By symmetry, we may assume that $u$, $v$, and $w$ appear in this order within the cyclic order on $C$.
If $w=v^+$, then the cycle $v^+\ldots u^-uxvv^-\ldots u^+v^+$ contradicts our assumption.
Hence, $w\not=v^+$.
If $w^-$ and $w^+$ are adjacent, 
then the cycle $w^-w^+\ldots u^-wvxuu^+\ldots v^-v^+\ldots w^-$ contradicts our assumption.
Hence, $w^-$ and $w^+$ are not adjacent.
If $w^+$ and $u$ are adjacent, 
then the cycle $uw^+\ldots u^-u^+\ldots v^-v^+\ldots w^-wvxu$ contradicts our assumption.
Hence, $w^+$ and $u$ are not adjacent,
which implies that $w^+\not=u^-$.
If $w^-$ and $u$ are adjacent, 
then the cycle $uw^- \ldots v^+v^-\ldots u^+u^-\ldots w^+wvxu$ contradicts our assumption.
Hence, $w^-$ and $u$ are not adjacent,
which implies that $w^-\not=u^+$.
If $w^+$ and $v$ are adjacent, 
then the cycle $uxvw^+\ldots u^-u^+\ldots v^-v^+\ldots wu$ contradicts our assumption.
Hence, $w^+$ and $v$ are not adjacent.
If $w^-$ and $v$ are adjacent, 
then the cycle $uxvw^-\ldots v^+v^-\ldots u^+u^-\ldots w^+wu$ contradicts our assumption.
Hence, $w^-$ and $v$ are not adjacent.
If $u'\in \{ u^-,u^+\}$ and $w'\in \{ w^-,w^+\}$, 
then $G[\{ u,v,w,u',w'\}]$ is not $(K_1+(K_1\cup P_3))$,
which implies that $u'$ and $w'$ are adjacent.
Hence, every vertex in $\{ u^-,u^+\}$ is adjacent to every vertex in $\{ w^-,w^+\}$.
Now, $G[\{ u,u^-,u^+,w^-,w^+\}]$ is $K_{1,1,3}$, which is a contradiction,
and completes the proof. $\Box$

\medskip

\noindent As we have seen in the previous results, the graph $K_{1,1,3}$ plays a special role.
Excluding this single graph, we obtain the following result.

\begin{theorem}\label{theorem2}
Every connected locally connected $\{ K_1+P_4,K_{1,4}, K_2+(K_1\cup K_2)\}$-free graph of order at least $3$ 
that is distinct from $K_{1,1,3}$ is fully cycle extendable.
\end{theorem}
{\it Proof:}
Let $G$ be a connected locally connected $\{ K_1+P_4,K_{1,4}, K_2+(K_1\cup K_2)\}$-free graph of order at least $3$.
As before,
the minimum degree of $G$ is at least $2$, and every vertex of $G$ lies on a triangle.
For a contradiction, we may assume that $G$ is distinct from $K_{1,1,3}$,
and that $C$ is a cycle in $G$ of order less than $n(G)$
such that $G$ contains no cycle $C'$ of order $n(C)+1$ with $V(C)\subseteq V(C')$.
We fix a cyclic order on $C$.
For every vertex $u$ on $C$, 
let $u^-$ and $u^+$ be the predecessor and successor of $u$ on $C$ within the cyclic order.
Since $G$ is connected, some vertex $u$ on $C$ has a neighbor $x$ in $V(G)\setminus V(C)$.
Our assumption implies that $x$ is adjacent to neither $u^-$ nor $u^+$.
Since $G[N_G(u)]$ is connected and $P_4$-free, 
$G[N_G(u)]$ contains a path $xvu^-$.
Our assumption implies that $v$ lies on $C$.

First, we assume that $u^-$ is adjacent to $u^+$.
Now, 
if $v$ is not adjacent to $u^+$, then $G[\{ u,u^-,u^+,v,x\}]$ is $K_1+P_4$, and,
if $v$ is adjacent to $u^+$, then $G[\{ u,u^-,u^+,v,x\}]$ is $K_2+(K_1\cup K_2)$, 
which is a contradiction.
Hence, $u^-$ is not adjacent to $u^+$.
In fact, by symmetry, for every vertex $\tilde{u}$ on $C$ that has a neighbour in $V(G)\setminus V(C)$,
we obtain that $\tilde{u}^-$ is not adjacent to $\tilde{u}^+$.
In particular, $v^-$ is not adjacent to $v^+$.
If $v^-$ is adjacent to $u^-$, 
then the cycle $uxvv^+\ldots u^-v^-\ldots u^+u$ contradicts our assumption.
Hence, $v^-$ and $u^-$ are not adjacent.

Next, we assume that $v^+\not=u^-$.
Since $G[\{ u^-,x,v,v^-,v^+\}]$ is not $K_{1,4}$,
$u^-$ and $v^+$ are adjacent.
If $u$ is not adjacent to $v^+$, then $G[\{ u,u^-,v,v^+,x\}]$ is $K_1+P_4$, and,
if $u$ is adjacent to $v^+$, then $G[\{ u,u^-,v,v^+,x\}]$ is $K_2+(K_1\cup K_2)$, 
which is a contradiction.
Hence, $v^+=u^-$.

Next, we assume that $v$ is adjacent to $u^+$.
By symmetry between $u^-$ and $u^+$, we obtain $v^-=u^+$, 
that is, $C$ has order $4$.
Since $G$ is connected and distinct from $K_{1,1,3}$,
there is a vertex $y\not\in \{ u,v,u^-,u^+,x\}$ with a neighbor in $\{ u,v,u^-,u^+,x\}$.
If $y$ is adjacent to $u$ or $v$, 
then, by symmetry, and since $G$ is $K_{1,4}$-free,
we obtain that $N_G(y)\cap \{ u,v,u^-,u^+,x\}=\{ u,v,x\}$,
and $G[\{ u,v,u^+,x,y\}]$ is $K_2+(K_1\cup K_2)$,
which is a contradiction.
Hence, $y$ is not adjacent to $u$ or $v$.
If $y$ is adjacent to $u^+$ or $u^-$, 
then we obtain, by symmetry, that $u^+$ is neighbor of $u^-$,
which is a contradiction.
Hence, the only neighbor of $y$ in $\{ u,v,u^-,u^+,x\}$ is $x$.
Since $G$ is locally connected and $(K_1+P_4)$-free, 
this implies the existence of a vertex $z\not\in \{ u,v,u^-,u^+,x\}$
that is adjacent to $x$, $y$, and $u$.
Since $z\not\in \{ u,v,u^-,u^+,x\}$, and $z$ is adjacent to $u$,
we obtain a similar contradiction as above.
Hence, $v$ is not adjacent to $u^+$.

By symmetry, there is a path $xwu^+$ in $G[N_G(u)]$
such that $w^-=u^+$, and $w$ is not adjacent to $u^-$.
If $v$ and $w$ are adjacent, 
then $G[\{ u,u^-,u^+,v,w\}]$ is $K_1+P_4$, and, 
if $v$ and $w$ are not adjacent, 
then $G[\{ u,u^-,v,x,w\}]$ is $K_1+P_4$, 
which is a contradiction, 
and completes the proof.
$\Box$

\medskip

\noindent We proceed to our results on locally Ore and locally Dirac graphs. 
As we will see below, 
the neighborhood condition used in the next result is weaker than being locally Ore or locally Dirac.

\begin{theorem}\label{theorem6}
If $G$ is a connected graph of order at least $3$ such that 
\begin{eqnarray}\label{e1}
|N_G(u)\cap N_G(v)\cap N_G(w)|
>|N_G(u)\setminus (N_G[v]\cup N_G[w])|
\end{eqnarray}
for every induced path $vuw$ of order $3$ in $G$,
then $G$ is fully cycle extendable.
\end{theorem}
{\it Proof:} Let $G$ be as in the statement.
Clearly, (\ref{e1}) implies that every vertex of $G$ of degree at least $2$ lies on a triangle. If $v$ is a vertex of degree $1$ in $G$, and $u$ is the unique neighbor of $v$, then the connectivity and $n(G)\geq 3$ imply that $u$ has another neighbor $w$,
and $vuw$ is an induced path of order $3$ that violates (\ref{e1}). Therefore, $G$ has minimum degree at least $2$, and every vertex of $G$ lies on a triangle.
For a contradiction, we may assume that $C$ is a cycle in $G$ of length less than $n(G)$ such that $G$ does not contain a cycle $C'$ of length $n(C)+1$ with $V(C)\subseteq V(C')$.
We fix a cyclic order on $C$.
For every vertex $u$ on $C$, 
let $u^+$ be the successor of $u$ on $C$ within the cyclic order.

Since $G$ is connected, some vertex $z$ in $V(G)\setminus V(C)$ has a neighbor on $C$.

Let $u$ be some neighbor of $z$ on $C$.
By our assumption, $z$ is not adjacent to $u^+$,
and $N_G(u)\cap N_G(u^+)\cap N(z)\subseteq V(C)$.
Let
\begin{eqnarray*}
A(u) &=& \Big\{v\in N_G(u)\cap N_G(u^+)\cap N(z): v^+\not\in N_G(u)\Big\},\\
B(u) &=& \Big\{v\in N_G(u)\cap N_G(u^+)\cap N(z): v^+\in N_G(u)\Big\},\mbox{ and}\\
C(u) &=& \Big\{v\in (N_G(u)\cap N(z)\cap V(C))\setminus N_G(u^+): v^+\in N_G(u)\Big\}.
\end{eqnarray*}
Note that $A(u)$, $B(u)$, and $C(u)$
are disjoint subsets of $V(C)$.

If $x$ and $y$ are distinct neighbors of $z$ on $C$, 
then our assumption implies that $x^+$ and $y^+$ are not adjacent.
Therefore, if $v\in B(u)\cup C(u)$, 
then $v^+$ is not adjacent to $z$ or $u^+$,
and hence
$v^+\in N_G(u)\setminus (N_G[u^+]\cup N_G[z])$,
which implies
$|N_G(u)\setminus (N_G[u^+]\cup N_G[z])|\geq 
|B(u)\cup C(u)|.$
Now, 
\begin{eqnarray*}
|A(u)|+|B(u)| & = & |N_G(u)\cap N_G(u^+)\cap N(z)|\\
& \stackrel{(\ref{e1})}{>} & 
|N_G(u)\setminus (N_G[u^+]\cup N_G[z])|\\
&\geq & |B(u)\cup C(u)|\\
& = & |B(u)|+|C(u)|,
\end{eqnarray*}
which implies 
\begin{eqnarray}\label{e2}
|A(u)| &>&|C(u)|
\end{eqnarray}
for every neighbor $u$ of $z$ on $C$.

Now, we define a sequence 
$u_1,u_2,u_3,\ldots$
of not necessarily distinct neighbors of $z$ on $C$.
Furthermore, for every neighbor $x$ of $z$ on $C$, 
we define two sets 
$A_k(x)\subseteq A(x)$
and 
$C_k(x)\subseteq C(x)$ 
for every positive integer $k$ 
for which $u_k$ is defined
in such a way that 
$$A_1(x)\subseteq A_2(x)\subseteq A_3(x)\ldots 
\subseteq A(x)$$
and
$$C_1(x)\subseteq C_2(x)\subseteq C_3(x)\ldots 
\subseteq C(x).$$
Let $u_1$ be any neighbor of $z$ on $C$,
and let $A_1(x)=C_1(x)=\emptyset$ for every neighbor $x$ of $z$ on $C$.
Now, we assume that 
the vertices $u_1,\ldots,u_k$
as well as 
the sets $A_k(x)$
and $C_k(x)$
for the neighbors $x$ of $z$ on $C$
have already been defined for some positive integer $k$.

If $A(u_k)\setminus A_k(u_k)\not=\emptyset$, 
then let $u_{k+1}\in A(u_k)\setminus A_k(u_k)$.
Note that, by definition, 
$u_{k+1}\in A(u_k)$ implies that $u_k\in C(u_{k+1})$.
Let 
$$
A_{k+1}(x)=
\left\{
\begin{array}{ll}
A_k(u_k)\cup \{ u_{k+1}\} &, x=u_k\mbox{ and}\\
A_k(x) &,\mbox{ otherwise},
\end{array}
\right.
$$
and
$$C_{k+1}(x)=
\left\{
\begin{array}{ll}
C_k(u_{k+1})\cup \{ u_k\} &, x=u_{k+1}\mbox{ and}\\
C_k(x) &,\mbox{ otherwise}.
\end{array}
\right.
$$
By the choice of $u_{k+1}$, 
we have $u_{k+1}\not\in A_k(u_k)$, 
and hence $|A_{k+1}(u_k)|=|A_k(u_k)|+1$.
If $u_k\in C_k(u_{k+1})$, 
then the above definitions imply 
the existence of some integer $i$ with $1\leq i<k$ 
such that $u_i=u_k$ and $u_{i+1}=u_{k+1}$.
Now, we obtain the contradiction
$$u_{k+1}=u_{i+1}\in 
A_i(u_i)\cup \{ u_{i+1}\}=
A_{i+1}(u_i)=A_{i+1}(u_k)\subseteq A_k(u_i).$$
This implies $u_k\not\in C_k(u_{k+1})$, 
and hence $|C_{k+1}(u_{k+1})|=|C_k(u_{k+1})|+1$.
By a simple inductive argument, 
we obtain that, 
for every positive integer $k$
for which $u_k$ is defined,
\begin{itemize}
\item $$\sum_{x\in N_G(z)\cap V(C)}\left(|A_k(x)|+|C_k(x)|\right)=2(k-1),$$
\item if $u_k=u_1$, then
\begin{eqnarray*}
|A_k(x)| &=& |C_k(x)|\mbox{ for $x\in N_G(z)\cap V(C)$},
\end{eqnarray*}
and, 
\item 
if $u_k\not=u_1$, then 
\begin{eqnarray*}
|A_k(u_1)| &=& |C_k(u_1)|+1,\\
|A_k(u_k)| &=& |C_k(u_k)|-1,\mbox{ and}\\
|A_k(x)| &=& |C_k(x)|\mbox{ for $x\in (N_G(z)\cap V(C))\setminus \{ u_1,u_k\}$}.
\end{eqnarray*}
\end{itemize}
Since $|A_k(u_k)|\leq |C_k(u_k)|$ holds in every case,
we obtain
\begin{eqnarray*}
|A(u_k)\setminus A_k(u_k)| 
&=& |A(u_k)|-|A_k(u_k)|\\
&\geq & |A(u_k)|-|C_k(u_k)|\\
&\geq & |A(u_k)|-|C(u_k)|\\
& \stackrel{(\ref{e2})}{>} & 0
\end{eqnarray*}
for every positive integer $k$
for which $u_k$ is defined.
This implies that $u_1,u_2,u_3,\ldots$ 
is actually an infinite sequence,
that is, $u_k$ is defined for every positive integer $k$.

Since
$$
\sum_{x\in N_G(z)\cap V(C)}\left(|A_k(x)|+|C_k(x)|\right)
\leq \sum_{x\in N_G(z)\cap V(C)}|A(x)\cup C(x)|
\leq \sum_{x\in V(C)}n(C)= n(C)^2,$$
we obtain a contradiction for $k>n(C)^2$,
which completes the proof. $\Box$

\medskip

\noindent Note that 
Corollary \ref{corollary1}
extends the main result, 
Theorem 3.3, 
of Kubicka et al. \cite{kubicka}.

\begin{corollary}\label{corollary1}
Every connected locally Ore graph $G$ of order at least $3$ is fully cycle extendable.
\end{corollary}
{\it Proof:} Let $G$ be as in the statement.
Let $vuw$ be an induced path of order $3$ in $G$.
By inclusion-exclusion,
we obtain 
\begin{eqnarray*}
d_G(u) & = & 
|\{ v,w\}|+|(N_G(v)\cup N_G(w))\cap N_G(u)|+|N_G(u)\setminus (N_G[v]\cup N_G[w])|\\
&=& 2+\Big(|N_G(u)\cap N_G(v)|+|N_G(u)\cap N_G(w)|-|N_G(u)\cap N_G(v)\cap N_G(w)|\Big)\\
&& +|N_G(u)\setminus (N_G[v]\cup N_G[w])|.
\end{eqnarray*}
Since $G$ is locally Ore, we have
$|N_G(u)\cap N_G(v)|+|N_G(u)\cap N_G(w)|\geq d_G(u)$,
which implies
$|N_G(u)\cap N_G(v)\cap N_G(w)|
\geq 2+|N_G(u)\setminus (N_G[v]\cup N_G[w])|$,
and hence (\ref{e1}). 
Now the desired results follows from Theorem \ref{theorem6}. $\Box$

\medskip

\noindent Corollary \ref{corollary1} immediately implies the following.

\begin{corollary}\label{corollary2}
Every connected locally Dirac graph $G$ of order at least $3$ is fully cycle extendable.
\end{corollary}

\end{document}